\documentclass[preprint,12pt,authoryear]{elsarticle}

\usepackage{amssymb}
\usepackage{amsmath}

\usepackage{natbib}
\usepackage{enumerate}
\usepackage{bbm}
\usepackage{bm}
\usepackage{latexsym}

\usepackage{color}
\usepackage{threeparttable}
\usepackage{url}

\newtheorem{theorem}{Theorem}
\newtheorem{lemma}{Lemma}

\journal{arXiv}

\begin{document}

\begin{frontmatter}



\title{A note on blinded continuous monitoring for continuous outcomes} 


%
\author[label]{Long-Hao Xu} 
\ead{long-hao.xu@med.uni-goettingen.de}
\author[label]{Tim Friede} 

\affiliation[label]{organization={Department of Medical Statistics, University Medical Center Göttingen},
            addressline={Humboldtallee~32},
            city={Göttingen},
            postcode={37073},
            country={Germany}}


\begin{abstract}
Continuous monitoring is becoming more popular due to its significant benefits, including reducing sample sizes and reaching earlier conclusions. In general, it involves monitoring nuisance parameters (e.g., the variance of outcomes) until a specific condition is satisfied. The blinded method, which does not require revealing group assignments, was recommended because it maintains the integrity of the experiment and mitigates potential bias. Although \cite{FriedeMiller:2012} investigated the characteristics of blinded continuous monitoring through simulation studies, its theoretical properties are not fully explored. In this paper, we aim to fill this gap by presenting the asymptotic and finite-sample properties of the blinded continuous monitoring for continuous outcomes. Furthermore, we examine the impact of using blinded versus unblinded variance estimators in the context of continuous monitoring. Simulation results are also provided to evaluate finite-sample performance and to support the theoretical findings.
\end{abstract}

%


\end{frontmatter}



\section{Introduction} 

The behavior of the optimal sample size is commonly considered in testing the two-sample problem.
Assume the data $X_1,...,X_n$ are independent and identically distributed, following a normal distribution with mean $\mu_1$ and variance $\sigma^2$.
Similarly, assume the data $Y_1,...,Y_n$ are independent and identically distributed, following a normal distribution with mean $\mu_2$ and variance $\sigma^2$.
For simplicity, we assume that outcomes are observed in pairs $(X_i,Y_i)$, $i=1,...,n$. 
The related one-sided hypothesis test is given by the null hypothesis $H_0: \mu_1-\mu_2=0$ versus the alternative hypothesis $H_1: \mu_1-\mu_2>0$.
This test is designed to have a one-sided significance level $\alpha$ under $H_0$, and a power $1-\beta$ for an assumed group difference $\delta_a$.

In the context of the two-sample problem, distinguishing between blinded and unblinded methods is essential. In general, blinded methods do not require revealing group assignments, with the aim of maintaining the integrity of the experiment and reducing potential conscious or unconscious bias.
In recent years, blinded methods have attracted growing attention in both practical applications and theoretical research, particularly within the context of clinical trials (e.g., \cite{Friedeetal:2019}, \cite{IgetaMatsui:2022}). An increasing number of researchers are examining how the use of blinded methods affects final results, and are comparing blinded and unblinded approaches when addressing the same problem (e.g., \cite{Graylingetal:2018}, \cite{TarimaFlournoy:2022}).

For convenience, let $\bm{Z}_n=(Z_1,...,Z_{2n})$ denote the blinded sample, where the group assignments are unknown.
In this setup, half of the elements in $\bm{Z}_n$ correspond to $(X_1, \ldots, X_n)$ and the other half to $(Y_1, \ldots, Y_n)$.
Define $\bar{X}_n=(\sum_{i=1}^nX_i)/n$ and $\bar{Y}_n=(\sum_{i=1}^nY_i)/n$.
As a consequence, the blinded sample mean is denoted as $\bar{Z}_n=\sum_{i=1}^{2n}Z_i/(2n)=\sum_{i=1}^n(X_i+Y_i)/(2n)=(\bar{X}_n+\bar{Y}_n)/2$.
While $\bar{X}_n$ and $\bar{Y}_n$ cannot be explicitly computed in the blinded setting, they can still be derived theoretically.
Based on this, the blinded variance estimator is defined as
\begin{align*}
  \hat{\sigma}_{n,blind}^2 &= \frac{1}{2n-1}\sum_{i=1}^{2n}(Z_i-\bar{Z}_n)^2 \\
  &= \frac{1}{2n-1}\left(\sum_{i=1}^n(X_i-\bar{Z}_n)^2+\sum_{j=1}^n(Y_j-\bar{Z}_n)^2\right) \\
  &= \frac{1}{2n-1}\left(\sum_{i=1}^n\left(X_i-\bar{X}_n+\frac{\bar{X}_n-\bar{Y}_n}{2}\right)^2+\sum_{j=1}^n\left(Y_j-\bar{Y}_n+\frac{\bar{Y}_n-\bar{X}_n}{2}\right)^2\right) \\
  &= \frac{1}{2n-1}\left(\sum_{i=1}^n(X_i-\bar{X}_n)^2+\sum_{j=1}^n(Y_j-\bar{Y}_n)^2+\frac{n(\bar{X}_n-\bar{Y}_n)^2}{2}\right).
\end{align*}
The quantity $\frac{1}{2n-1}(\sum_{i=1}^n(X_i-\bar{X}_n)^2+\sum_{j=1}^n(Y_j-\bar{Y}_n)^2+\frac{n(\bar{X}_n-\bar{Y}_n)^2}{2})$ is an idealized mathematical expression that cannot be computed in the blinded setting. Instead, we use $\frac{1}{2n-1}\sum_{i=1}^{2n}(Z_i-\bar{Z}_n)^2$ to calculate $\hat{\sigma}_{n,blind}^2$.
For comparison, the unblinded variance estimator is defined as
\begin{align*}
  \hat{\sigma}_{n,unblind}^2 &= \frac{1}{2n-2}\left(\sum_{i=1}^n(X_i-\bar{X}_n)^2+\sum_{j=1}^n(Y_j-\bar{Y}_n)^2\right) \\
  &= \frac{1}{2}\left(\frac{1}{n-1}\sum_{i=1}^n(X_i-\bar{X}_n)^2+\frac{1}{n-1}\sum_{j=1}^n(Y_j-\bar{Y}_n)^2\right).
\end{align*}

In contrast to the unblinded variance estimator $\hat{\sigma}_{n,unblind}^2$, the blinded variance estimator $\hat{\sigma}_{n,blind}^2$ can be computed without knowledge of the group assignments.
Note that $\hat{\sigma}_{n,unblind}^2$ is unbiased and consistent, while $\hat{\sigma}_{n,blind}^2$ is biased and not consistent.
For more details on the statistical properties of $\hat{\sigma}_{n,unblind}^2$ and $\hat{\sigma}_{n,blind}^2$, readers are referred to \cite{FriedeKieser:2013}.

If the true variance $\sigma^2$ were known, the optimal sample size required to achieve a significance level $\alpha$ under $H_0$ and power $1-\beta$ at the assumed group difference $\delta_a$ would be $n_{req}=v\sigma^2$ per group, where $v=2(\Phi^{-1}(1-\alpha)+\Phi^{-1}(1-\beta))^2/\delta_a^2$. Detailed derivations can be found, for example, in Chapter 3 of \cite{Kieser:2020}.
This quantity $v$ essentially reflects the required sample size per unit of variance, and it tells us how many samples are needed per unit variance to achieve the desired power and significance level.
In most situations, especially in clinical trials, the true variance $\sigma^2$ is unknown and must be estimated during the study.

Continuous monitoring (also known as a form of sequential analysis) has been extensively studied since the pioneering work of \cite{Anscombe:1953} and \cite{Robbins:1959}.
Subsequent contributions by, for example, \cite{Woodroofe:1977}, \cite{Siegmund:1985}, \cite{Ghoshetal:1997}, \cite{Lai:2001}, and \cite{MukhopadhyaySilva:2008} have provided a theoretical foundation in this area. 
These studies primarily focused on unblinded procedures, which typically rely on known group assignments to estimate nuisance parameters such as the variance.
In contexts where unblinding is not permitted, such as during interim analyses in clinical trials, blinded methods become essential. Addressing this gap,
\cite{FriedeMiller:2012} proposed the following sequential procedure, referred to as blinded continuous monitoring, and they defined the number of observations per group as
\begin{align}\label{N-process}
N_b=\min\{ n=n_1,n_1+1,... | \hat{\sigma}^2_{n,blind}\leq n/v \}
\end{align}
where $\hat{\sigma}^2_{n,blind}$ is the blinded variance estimator for the current observed samples, $n_1$ is the initial sample size, and $v=2(\Phi^{-1}(1-\alpha)+\Phi^{-1}(1-\beta))^2/\delta_a^2$. 
Similarly, when using the unblinded variance estimator $\hat{\sigma}_{n,unblind}^2$, the number of observations per group is defined as
\begin{align}\label{N-process-un}
N_u=\min\{ n=n_1,n_1+1,... | \hat{\sigma}^2_{n,unblind}\leq n/v \}.
\end{align}

As defined in (\ref{N-process}), blinded continuous monitoring can generally be interpreted as sequentially evaluating the blinded variance estimator $\hat{\sigma}^2_{n,blind}$ at each sample size $n$ until a predefined stopping criterion is met. The definition and interpretation of $N_u$ are analogous.
Generally, continuous monitoring aims to reduce the total sample size by evaluating the nuisance parameter(s) until predefined conditions are met.
In this framework, the sample size is not fixed in advance, and a decision may be reached either earlier or later than in a fixed-sample design where the true variance $\sigma^2$ is known.

Although the unblinded variance estimator $\hat{\sigma}^2_{n,unblind}$ is widely used in many fields, regulatory agencies such as the U.S. Food and Drug Administration and the European Medicines Agency recommend the use of blinded methods during interim analyses in clinical trials whenever feasible.
In such cases, the blinded variance estimator $\hat{\sigma}^2_{n,blind}$ is a necessary alternative.
However, theoretical results for the blinded continuous monitoring remain largely unexplored. 
To the best of our knowledge, the blinded method had not been studied in the context of sequential analysis until the work of \cite{FriedeMiller:2012}, which primarily focused on numerical results.
\cite{Liu:1997} considered a similar question, but he only discussed the one-sample problem and he did not mention the blinded method.
\cite{HuMukhopadhyay:2019} provided general conditions that general variance estimators should satisfy to ensure desirable properties.
However, the blinded method was not considered in their framework and, in fact, the blinded variance estimator $\hat{\sigma}^2_{n,blind}$ does not meet their stated conditions. 
Given the growing importance of blinded methods, it is essential to investigate their statistical properties when applied in the context of continuous monitoring.

In this paper, we explore the application of the blinded method in clinical trials, with a particular focus on how continuous monitoring is affected by its use.
First, we provide a theoretical proof of the results in \citet{FriedeMiller:2012}, giving the asymptotic and finite-sample properties of blinded continuous monitoring (\ref{N-process}).
We then investigate how different limiting situations of $v$ or $\sigma$ affect these asymptotic properties.
Our findings show that the blinded method affects the asymptotic properties of continuous monitoring, yielding different results as $v\rightarrow\infty$ compared to $\sigma\rightarrow\infty$.
In contrast, no such difference is observed for the unblinded continuous monitoring (\ref{N-process-un}).

The remainder of the paper is organized as follows. Section \ref{sec-2} examines the asymptotic and finite-sample properties of the blinded continuous monitoring (\ref{N-process}).
Section \ref{sec-3} illustrates the simulation results.
Remarks and discussions are presented in Section \ref{sec-4}. All proofs are provided in the supplementary material.


\section{Main results}\label{sec-2} 

In this section, we introduce several properties of the blinded continuous monitoring (\ref{N-process}).
\begin{theorem}\label{th-1}
If the initial sample size $n_1\geq2$, $0<\sigma^2<\infty$, and $0<v<\infty$, then we have (i) $N_b$ is well-defined and $\mathbb{P}(N_b<\infty)=1$, that is, the stochastic process $N_b$ terminates with probability 1.

(ii)
\begin{align*}
\mathbb{E}(N_b) \leq n_1 + v\sigma^2+\frac{v(\mu_1-\mu_2)^2}{4} = n_1 + n_{req}\left(1+\frac{(\mu_1-\mu_2)^2}{4\sigma^2}\right).
\end{align*}

(iii) 
\begin{align*}
\mathbb{E}(N_b^2) \leq \left( n_1 + v\sigma^2+\frac{v(\mu_1-\mu_2)^2}{4} \right)^2 = \left( n_1 + n_{req}\left(1+\frac{(\mu_1-\mu_2)^2}{4\sigma^2}\right) \right)^2.
\end{align*}
\end{theorem}

Note that the upper bound of $\mathbb{E}(N_b)$ and $\mathbb{E}(N_b^2)$ depend on $v$, the true variance $\sigma^2$ and the true group difference $\mu_1-\mu_2$.
The term $(\mu_1-\mu_2)^2/(4\sigma^2)$ shown in Theorem \ref{th-1} (ii) and (iii) can be interpreted as the cost of using blinded continuous monitoring, and it is typically less than 1 in most clinical trials (cf. \cite{KieserFriede:2003}).
From Theorem \ref{th-1} (iii), an upper bound of $Var(N_b)$ can also be derived, showing that the variability in the blinded continuous monitoring remains controlled.

The following asymptotic property of the blinded continuous monitoring (\ref{N-process}) highlights its conservative nature. It is unlikely to underestimate $n_{req}$, as $\sigma\rightarrow\infty$ or $v\rightarrow\infty$.

\begin{theorem}\label{th-P(N<epsilona)}
For every $\varepsilon$ in $(0,1)$, if the initial sample size $n_1\geq2$, then $\mathbb{P}(N_b\leq \varepsilon n_{req})  = \mathcal{O}((n_{req})^{-(n_1-1)})$ as $\sigma\rightarrow\infty$ or $v\rightarrow\infty$.
\end{theorem}

For other asymptotic properties, we consider the cases where $v\rightarrow\infty$ and $\sigma\rightarrow\infty$, respectively.
When unblinded continuous monitoring (\ref{N-process-un}) is considered, the asymptotic properties of (\ref{N-process-un}) will align with Theorem \ref{th-2}, regardless of whether $v\rightarrow\infty$ and $\sigma\rightarrow\infty$ (cf. Section 13.3 in \cite{MukhopadhyaySilva:2008}).
In contrast, the asymptotic results for the blinded continuous monitoring (\ref{N-process}) are different between the two limiting scenarios $v\rightarrow\infty$ and $\sigma\rightarrow\infty$ due to the influence of the blinded method.
First, we consider the case where $\sigma\rightarrow\infty$ indicating that the variance of both $X_i$ and $Y_i$ approaches infinity.
In this case, their observed values become highly dispersed around their respective means, making statistical inference more challenging.
Second, we consider the case where $v\rightarrow\infty$.
Given a fixed Type I error rate $\alpha$ and power $1-\beta$ at an assumed group difference $\delta_a$, the condition $v\rightarrow\infty$ corresponds to $\delta_a\rightarrow0$, 
meaning that the alternative hypothesis approaches the null hypothesis, i.e., the group difference becomes arbitrarily small.
In this case, the power $1-\beta$ is achieved at the assumed group difference $\delta_a$.
The values of $v$, $\mu_1$, and $\mu_2$ are fixed and finite in Theorem \ref{th-2}, while the values of $\sigma$, $\mu_1$, and $\mu_2$ are fixed and finite in Theorem \ref{th-3}.

\begin{theorem}\label{th-2}

(i)
$N_b$ is well-defined and non-decreasing as a function of $\sigma$. Furthermore, we have $N_b(\sigma)\rightarrow\infty$ almost surely as $\sigma\rightarrow\infty$, and $\mathbb{E}(N_b)\rightarrow\infty$ as $\sigma\rightarrow\infty$.

(ii)
\begin{align*}
\lim_{\sigma\rightarrow\infty}\frac{N_b}{n_{req}} = 1 \quad almost \ surely.
\end{align*}

(iii)
\begin{align*}
\lim_{\sigma\rightarrow\infty} \mathbb{E}\left(\frac{N_b}{n_{req}}\right) = 1.
\end{align*}

(iv)
\begin{align*}
\lim_{\sigma\rightarrow\infty}\mathbb{E}\left(\frac{N_b^2}{n_{req}^2}\right) = 1.
\end{align*}

(v) 
As $\sigma\rightarrow\infty$, we have
\begin{align*}
\frac{N_b-n_{req}}{\sqrt{n_{req}}} \overset{d}{\rightarrow} N(0,1).
\end{align*}
\end{theorem}

Theorem \ref{th-2} (i) shows that both the random sample size $N_b$ and the expected sample size $\mathbb{E}(N_b)$ diverge to infinity as $\sigma\rightarrow\infty$, which means that larger sample sizes are required to maintain statistical precision in high-variance settings.
Theorem \ref{th-2} (ii) is referred to as asymptotic optimality.
This property implies that $N_b$ will almost surely be closed to $n_{req}$ with the increasing of $\sigma$. 
Theorem \ref{th-2} (iii) is referred to as asymptotic efficiency.
This means that the expected sample size $\mathbb{E}(N_b)$ will be closed to $n_{req}$ with the increasing of $\sigma$ and the procedure remains asymptotic efficient in terms of expected sample size compared to $n_{req}$ in the fixed-sample design.
Theorem \ref{th-2} (iv) provides insights into the second-moment properties, which is directly related to the variance of $N_b$ as $\sigma\rightarrow\infty$.
Theorem \ref{th-2} (v) establishes the asymptotic normality of $N_b$ and indicates that $(N_b-n_{req})/\sqrt{n_{req}}$ follows an approximately normal distribution when $\sigma$ is sufficiently large.

\begin{theorem}\label{th-3}

(i)
$N_b$ is well-defined and non-decreasing as a function of $v$. Furthermore, we have $N_b(v)\rightarrow\infty$ almost surely as $v\rightarrow\infty$, and $\mathbb{E}(N_b)\rightarrow\infty$ as $v\rightarrow\infty$.

(ii)
\begin{align*}
\lim_{v\rightarrow\infty}\frac{N_b}{n_{req}} = 1+\frac{(\mu_1-\mu_2)^2}{4\sigma^2} \quad almost \ surely.
\end{align*}

(iii)
\begin{align*}
\lim_{v\rightarrow\infty} \mathbb{E}\left(\frac{N_b}{n_{req}}\right) = 1 + \frac{(\mu_1-\mu_2)^2}{4\sigma^2}.
\end{align*}

(iv)
\begin{align*}
\lim_{v\rightarrow\infty}\mathbb{E}\left(\frac{N_b^2}{n_{req}^2}\right) = \left(1+\frac{(\mu_1-\mu_2)^2}{4\sigma^2}\right)^2.
\end{align*}

(v) 
As $v\rightarrow\infty$, we have
\begin{align*}
\frac{N_b-n_{req}(1+(\mu_1-\mu_2)^2/(4\sigma^2))}{\sqrt{n_{req}}} \overset{d}{\rightarrow} N\left(0,\frac{4\sigma^2+2(\mu_1-\mu_2)^2}{4\sigma^2+(\mu_1-\mu_2)^2}\right).
\end{align*}
\end{theorem}

Theorem \ref{th-3} (i) shows that both the random sample size $N_b$ and the expected sample size $\mathbb{E}(N_b)$ diverge to infinity as $v\rightarrow\infty$, indicating that a larger sample size is required in this limiting case.
Theorem \ref{th-3} (ii) states that the random sample size $N_b$ asymptotically exceeds $n_{req}$, when the true group difference $\delta=\mu_1-\mu_2\neq0$.
Similarly, this term $(\mu_1-\mu_2)^2/(4\sigma^2)$ can be interpreted as the cost of using blinded continuous monitoring. 
From Theorem \ref{th-3} (iii), we know that the expected sample size $\mathbb{E}(N_b)$ gets very close to $n_{req}(1+(\mu_1-\mu_2)^2/(4\sigma^2))$ as $v$ increases, indicating that the blinded continuous monitoring leads to a slight inflation in sample size compared to the fixed-sample design.
Theorem \ref{th-3} (iv) shows that there is one more additional term $(\mu_1-\mu_2)^2/(4\sigma^2)$ appearing even in the second-moment properties, which is directly related to the variance of $N_b$ as $v\rightarrow\infty$.
In addition, Theorem \ref{th-3} (v) establishes that the random sample size $N_b$ asymptotically follows a normal distribution.
Specifically, it shows that the normalized stopping variable $(N_b-n_{req}(1+(\mu_1-\mu_2)^2/4\sigma^2))/\sqrt{n_{req}}$ would have an approximate normal distribution when $v$ is sufficiently large.

\section{Simulation}\label{sec-3}

In this section, we present the simulation results and evaluate the performance of blinded continuous monitoring (\ref{N-process}) and unblinded continuous monitoring (\ref{N-process-un}) for finite values $v$ and $\sigma$.
To ensure a fair comparison, we choose the optimal sample size in the fixed-sample design as $n_{req}=10,50,100,500,1000$ for both blinded and unblinded monitoring approaches.
To examine the differences between the two limiting situations, we keep $v=1$ as fixed and vary $\sigma$ across $\sqrt{10}$, $5\sqrt{2}$, 10, $10\sqrt{5}$, and $10\sqrt{10}$ in Table \ref{table-1}, while we keep $\sigma=1$ as fixed and vary $v$ across 10, 50, 100, 500, and 1000 in Table \ref{table-2}.
In both tables, we choose $\mu_1=0,2,5$ and keep $\mu_2$ fixed.
All simulations were conducted using R computing environment, version 4.2.0 \citep{R:2022}, each with 10000 simulations per scenario. 
In addition, we use the same random seed in each scenario.

For convenience, we introduce several notations used in Table \ref{table-1} and Table \ref{table-2}.
Denote $\overline{N}_b$ as the sample mean of the final sample size implementing blinded continuous monitoring of the variance for all 10000 simulations. The standard error of $\overline{N}_b$ is denoted as  $s(\overline{N}_b)$.
Similarly, denote $\overline{N}_u$ as the sample mean of the final sample size implementing unblinded continuous monitoring of the variance for all 10000 simulations.
The standard error of $\overline{N}_u$ is denoted as $s(\overline{N}_u)$.

\begin{table}[ht]
\centering
\begin{threeparttable}[b]
\scriptsize
\caption{Simulations from $N(\mu_1,\sigma^2)$ and $N(\mu_2,\sigma^2)$ under 10000 runs implementing blinded and unblinded continuous monitoring of the variance ($v=1$ and $\mu_2=0$ are fixed, the initial sample size $n_1=10$)}
\label{table-1}
\begin{tabular}{cccccccccc}
  \hline
  \hline
  $\mu_1$ & $n_{req}$ & $\sigma$ & $\overline{N}_b$ & $s(\overline{N}_b)$  & $\overline{N}_b/n_{req}$ & Upper bound\tnote{1} & $\overline{N}_u$ & $s(\overline{N}_u)$  & $\overline{N}_u/n_{req}$ \\
  \hline
  ~ & 10 & $\sqrt{10}$ & 11.4037 & 1.8735 & 1.1404 & 11.75 & 11.2955 & 1.8238 & 1.1296 \\
  ~ & 50 & $5\sqrt{2}$ & 49.9064 & 7.3498 & 0.9981 & 51.75 & 49.7201 & 7.3623 & 0.9944 \\
  $\mu_1=1$ & 100 & 10 & 99.9400 & 10.1021 & 0.9994 & 101.75 & 99.5718 & 10.1720 & 0.9957 \\
  ~ & 500 & $10\sqrt{5}$ & 500.1991 & 22.3504 & 1.0004 & 501.75 & 499.9163 & 22.1464 & 0.9998 \\
  ~ & 1000 & $10\sqrt{10}$ & 1000.0603 & 31.4952 & 1.0001 & 1001.75 & 999.9493 & 31.5001 & 0.9999 \\
  \hline
  ~ & 10 & $\sqrt{10}$ & 11.9195 & 2.2187 & 1.1920 & 12.5 & 11.2955 & 1.8238 & 1.1296 \\
  ~ & 50 & $5\sqrt{2}$ & 50.7311 & 7.4199 & 1.0146 & 52.5 & 49.7201 & 7.3623 & 0.9944 \\
  $\mu_1=2$ & 100 & 10 & 100.6331 & 10.1656 & 1.0063 & 102.5 & 99.5718 & 10.1720 & 0.9957 \\
  ~ & 500 & $10\sqrt{5}$ & 500.9419 & 22.3626 & 1.0019 & 502.5 & 499.9163 & 22.1464 & 0.9998 \\
  ~ & 1000 & $10\sqrt{10}$ & 1000.8652 & 31.5596 & 1.0009 & 1002.5 & 999.9493 & 31.5001 & 0.9999 \\
  \hline
  ~ & 10 & $\sqrt{10}$ & 16.4330 & 3.7084 & 1.6433 & 17.75 & 11.2955 & 1.8238 & 1.1296 \\\
  ~ & 50 & $5\sqrt{2}$ & 56.0726 & 7.6538 & 1.1215 & 57.75 & 49.7201 & 7.3623 & 0.9944 \\
  $\mu_1=5$ & 100 & 10 & 106.0004 & 10.3054 & 1.0600 & 107.75 & 99.5718 & 10.1720 & 0.9957 \\
  ~ & 500 & $10\sqrt{5}$ & 506.1834 & 22.5512 & 1.0124 & 507.75 & 499.9163 & 22.1464 & 0.9998 \\
  ~ & 1000 & $10\sqrt{10}$ & 1006.0055 & 31.5900 & 1.0060 & 1007.75 & 999.9493 & 31.5001 & 0.9999 \\
  \hline
  \hline
\end{tabular}
   \begin{tablenotes}
     \item[1] The upper bound of $\mathbb{E}(N_b)$ in Theorem \ref{th-1} in blinded continuous monitoring
   \end{tablenotes}
  \end{threeparttable}
\end{table}

\begin{table}[ht]
\centering
\begin{threeparttable}[b]
\scriptsize
\caption{Simulations from $N(\mu_1,\sigma^2)$ and $N(\mu_2,\sigma^2)$ under 10000 runs implementing blinded and unblinded continuous monitoring of the variance ($\sigma=1$ and $\mu_2=0$ are fixed, the initial sample size $n_1=10$)}
\label{table-2}
\begin{tabular}{cccccccccc}
  \hline
  \hline
  $\mu_1$ & $n_{req}$ & $v$ & $\overline{N}_b$ & $s(\overline{N}_b)$  & $\overline{N}_b/n_{req}$ & Upper bound\tnote{1} & $\overline{N}_u$ & $s(\overline{N}_u)$  & $\overline{N}_u/n_{req}$ \\
  \hline
  ~ & 10 & $ 10 $ & 13.0250 & 2.7759 & 1.3025 & 14 & 11.2955 & 1.8238 & 1.1296  \\
  ~ & 50 & $ 50 $ & 62.4143 & 7.9760 & 1.2483 & 64 & 49.7201 & 7.3623 & 0.9944   \\
  $\mu_1=1$ & 100 & 100 & 124.9483 & 10.9901 & 1.2495 & 126.5 & 99.5718 & 10.1720 & 0.9957  \\
  ~ & 500 & 500 & 624.8957 & 24.2892 & 1.2498 & 626.5 & 499.9163 & 22.1464 & 0.9998   \\
  ~ & 1000 & 1000 & 1249.7079 & 34.0695 & 1.2497 & 1251.5 & 999.9493 & 31.5001 & 0.9999  \\
  \hline
  ~ & 10 & $ 10 $ & 20.2562 & 4.0676 & 2.0256 & 21.5 & 11.2955 & 1.8238 & 1.1296  \\
  ~ & 50 & $ 50 $ & 100.2653 & 8.7786 & 2.0053 & 101.5 & 49.7201 & 7.3623 & 0.9944   \\
  $\mu_1=2$ & 100 & 100 & 200.3337 & 12.2293 & 2.0033 & 201.5 & 99.5718 & 10.1720 & 0.9957  \\
  ~ & 500 & 500 & 1000.0871 & 27.1989 & 2.0002 & 1001.5 & 499.9163 & 22.1464 & 0.9998   \\
  ~ & 1000 & 1000 & 2000.5408 & 38.4516 & 2.0005 & 2001.5 & 999.9493 & 31.5001 & 0.9999  \\
  \hline
  ~ & 10 & $ 10 $ & 73.3414 & 4.3575 & 7.3341 & 74 & 11.2955 & 1.8238 & 1.1296  \\
  ~ & 50 & $ 50 $ & 363.3671 & 9.6653 & 7.2673 & 364 & 49.7201 & 7.3623 & 0.9944   \\
  $\mu_1=5$ & 100 & 100 & 725.7457 & 13.6278 & 7.2575 & 726.5 & 99.5718 & 10.1720 & 0.9957  \\
  ~ & 500 & 500 & 3626.0818 & 30.3324 & 7.2522 & 3626.5 & 499.9163 & 22.1464 & 0.9998   \\
  ~ & 1000 & 1000 & 7251.0777 & 42.9622 & 7.2511 & 7251.5 & 999.9493 & 31.5001 & 0.9999  \\
  \hline
  \hline
\end{tabular}
   \begin{tablenotes}
     \item[1] The upper bound of $\mathbb{E}(N_b)$ in Theorem \ref{th-1} in blinded continuous monitoring
   \end{tablenotes}
  \end{threeparttable}
\end{table}

From Table \ref{table-1}, the average estimated sample size $\overline{N}_b$ and $\overline{N}_u$ lie within a small neighborhood of $n_{req}$.
In addition, the ratios $\overline{N}_b/n_{req}$ and $\overline{N}_u/n_{req}$ are close to 1 with the increase of $\sigma$, confirming the asymptotic efficiency properties of $N_{b}$ and $N_{u}$.
There is a relative small gap between the average estimated sample size $\overline{N}_b$ and the upper bound of $\mathbb{E}(N_b)$.
From Table \ref{table-2}, the average estimated sample size $\overline{N}_b$ tends to be larger than $n_{req}$, while the average estimated sample size $\overline{N}_u$ lies within a small neighborhood of $n_{req}$.
In addition, as $v$ increase, the ratio $\overline{N}_b/n_{req}$ approaches the theoretical value $(\mu_1-\mu_2)^2/(4\sigma^2)$, while the ratio $\overline{N}_u/n_{req}$ remains close to 1.
Similarly, there is a small gap between the average estimated sample size $\overline{N}_b$ and the upper bound of $\mathbb{E}(N_b)$ in Table \ref{table-2}.
In addition, we observe that $s(\overline{N}_b)$ tends to be larger than $s(\overline{N}_u)$.
It means that the blinded continuous monitoring offers a slightly larger variability than the unblinded continuous monitoring if there is a large $v$.
As $\mu_1$ increases, a larger sample size is required in the blinded continuous monitoring.

An interesting phenomenon is that the unblinded continuous monitoring (\ref{N-process-un}) yields the same results for a given $n_{req}$, regardless of the specific values of $v$, $\sigma$, $\mu_1$, and $\mu_2$, when the same random seed is used. We provide a heuristic explanation as follows.
Based on the definition of the unblinded variance estimator, it can be deduced that 
it follows that
\begin{align*}
  v\hat{\sigma}_{n,unblind}^2 &= \frac{v\sigma^2}{2n-2}\left(\sum_{i=1}^n\left(\frac{X_i-\bar{X}_n}{\sigma}\right)^2+\sum_{j=1}^n\left(\frac{Y_j-\bar{Y}_n}{\sigma}\right)^2\right) \\
  &= \frac{n_{req}}{2n-2}\left(\sum_{i=1}^n\left(\frac{X_i-\bar{X}_n}{\sigma}\right)^2+\sum_{j=1}^n\left(\frac{Y_j-\bar{Y}_n}{\sigma}\right)^2\right).
\end{align*}
Because $X_i\sim N(\mu_1,\sigma^2)$, $i=1,...,n$, and $Y_j\sim N(\mu_2, \sigma^2)$, $j=1,...,n$, we know $(X_i-\bar{X}_n)/\sigma \sim N(0,1+\frac{1}{n})$ and $(Y_j-\bar{Y}_n)/\sigma \sim N(0,1+\frac{1}{n})$.
Notice that the normal distribution $N(0,1+\frac{1}{n})$ are independent of $v$, $\sigma$, $\mu_1$, and $\mu_2$.
Therefore, the unblinded continuous monitoring (\ref{N-process-un}) only depends on $n_{req}$.
This implies that the simulation results will remain the same as long as $n_{req}$ does not change, provided the random seed is fixed.
However, the blinded continuous monitoring (\ref{N-process}) does not have this property. 
This also highlights the fundamental difference between blinded and unblinded continuous monitoring. The quantity $n_{req}$ has a significant impact on unblinded continuous monitoring. However, in blinded continuous monitoring, the process is influenced not only by $n_{req}$, but also by the true values of $v$, $\sigma$, and the mean difference $\mu_1-\mu_2$.

\section{Concluding remarks}\label{sec-4} 

In this paper, we investigate the impact of the blinded method in continuous monitoring proposed by \cite{FriedeMiller:2012}.
We provide several theoretical results related to blinded continuous monitoring.
The approach developed in this paper can be generalized to broader classes of blinded estimators and nuisance parameters within the context of sequential analysis. 
Moreover, it may inform sequential interval estimation and prediction interval methods, such as those proposed by \cite{IsogaiFutschik:2008}, \cite{Uno:2013}, and \cite{Schachtetal:2025}, where precise control of estimation accuracy is essential.

\section*{Acknowledgments}

The authors thank Ole Schacht, Lara Vankelecom, Tom Loeys, Beatrijs Moerkerke, and Kelly Van Lancker (all from Ghent University) for discussions on the topic that helped to improve the manuscript.
The authors are grateful for funding by the DFG grant numbers FR 3070/4-1 and FR 3070/4-2.

%
%




\newpage

{\centering {\Large A note on blinded continuous monitoring for continuous outcomes}}

\section*{Supplementary material}

Before providing the proofs of theorems, we first introduce the following useful lemma.

\begin{lemma}\label{lemma-1} 
Let $X_1,X_2,...$ be i.i.d. Define $S_n=\sum_{i=1}^n X_i(n\geq 1)$. Let $N$ be a stopping time with $\mathbb{E}(N)<\infty$. Then, if $\mathbb{E}(|X_1|)<\infty$, $\mathbb{E}(S_N)=\mu\mathbb{E}(N)$, where $\mu=\mathbb{E}(X_1)$.
\end{lemma}
\noindent\textbf{Proof of Lemma \ref{lemma-1}.} This is also called Wald's first lemma. See the detail of the proof in \citet[Theorem 2.4.4]{Ghoshetal:1997}.

\noindent\textbf{Proof of Theorem \ref{th-1}.}

(i) 

Because $0<\sigma^2<\infty$, we have $0\leq(\mu_1-\mu_2)^2<\infty$.
From the strong law of large numbers, we have 
\begin{align*}
\hat{\sigma}^2_{n,blind} \rightarrow \sigma^2 + (\mu_1-\mu_2)^2/4 \quad almost \ surely \  as \ n\rightarrow\infty.
\end{align*}

Since 
\begin{align*}
\lim_{n\rightarrow\infty} \mathbb{P}(N_b> n) \leq \lim_{n\rightarrow\infty}\mathbb{P}(\hat{\sigma}_{n,blind}^2 > n/v),
\end{align*}
it follows that
\begin{align*}
\mathbb{P}(N_b<\infty)
= \lim_{n\rightarrow\infty} \mathbb{P}(N_b\leq n) = 1.
\end{align*}

Thus, this completes the proof.

(ii)

First assume $\mathbb{E}(N_b)<\infty$.
From the definition of $N_b$, we must have
\begin{align}\label{eq-22}
N_b-1 \leq (n_1-1)\mathbb{I}_{N_b=n_1} + v\hat{\sigma}^2_{N_b-1,blind}.
\end{align}

Since (\ref{eq-22}) and the definition of $\hat{\sigma}^2_{N_b-1,blind}$, it follows that
\begin{align}\label{eq-23}
& (2N_b-3)(N_b-1) \nonumber \\
&\leq (2N_b-3)(n_1-1)\mathbb{I}_{N_b=n_1} + v\left(\sum_{i=1}^{N_b-1}(X_i-\bar{Z}_{N_b-1})^2+\sum_{j=1}^{N_b-1}(Y_j-\bar{Z}_{N_b-1})^2\right) \nonumber \\
&=(2n_1-3)(n_1-1)\mathbb{I}_{N_b=n_1} + v\left(\sum_{i=1}^{N_b-1}(X_i-\bar{Z}_{N_b-1})^2+\sum_{j=1}^{N_b-1}(Y_j-\bar{Z}_{N_b-1})^2\right) \nonumber \\
&=(2n_1-3)(N_b-1)\mathbb{I}_{N_b=n_1} + v\left(\sum_{i=1}^{N_b-1}(X_i-\bar{Z}_{N_b-1})^2+\sum_{j=1}^{N_b-1}(Y_j-\bar{Z}_{N_b-1})^2\right) \nonumber \\
&\leq (2n_1-3)(N_b-1) + v\left(\sum_{i=1}^{N_b-1}(X_i-\bar{Z}_{N_b-1})^2+\sum_{j=1}^{N_b-1}(Y_j-\bar{Z}_{N_b-1})^2\right) \nonumber \\
&\leq (2n_1-3)(N_b-1) + \nonumber \\
& v\left(\sum_{i=1}^{N_b-1}(X_i-\bar{Z}_{N_b-1}+\bar{Z}_{N_b-1}-(\mu_1+\mu_2)/2)^2+\sum_{j=1}^{N_b-1}(Y_j-\bar{Z}_{N_b-1}+\bar{Z}_{N_b-1}-(\mu_1+\mu_2)/2)^2\right) \nonumber \\
&= (2n_1-3)(N_b-1) + v\left(\sum_{i=1}^{N_b-1}(X_i-(\mu_1+\mu_2)/2)^2+\sum_{j=1}^{N_b-1}(Y_j-(\mu_1+\mu_2)/2)^2\right).
\end{align}

From Jensen's inequality, the left hand side of (\ref{eq-23}) leads to
\begin{align*}
(2\mathbb{E}(N_b)-3)(\mathbb{E}(N_b)-1) \leq \mathbb{E}((2N_b-3)(N_b-1)).
\end{align*}
From Lemma \ref{lemma-1}, the right hand side of (\ref{eq-23}) leads to
\begin{align*}
& (2n_1-3)\mathbb{E}(N_b-1)+v\mathbb{E}\left(\sum_{i=1}^{N_b-1}(X_i-(\mu_1+\mu_2)/2)^2+\sum_{j=1}^{N_b-1}(Y_j-(\mu_1+\mu_2)/2)^2\right) \\
&= (2n_1-3)(\mathbb{E}(N_b)-1) + 2v(\mathbb{E}(N_b)-1)\left(\sigma^2+\frac{(\mu_1-\mu_2)^2}{4}\right).
\end{align*}

Thus, we have
\begin{align*}
(2\mathbb{E}(N_b)-3)(\mathbb{E}(N_b)-1) \leq (2n_1-3)(\mathbb{E}(N_b)-1) + 2v(\mathbb{E}(N_b)-1)\left(\sigma^2+\frac{(\mu_1-\mu_2)^2}{4}\right),
\end{align*}
which completes the proof.

If we do not assume $\mathbb{E}(N_b)<\infty$, we define $N_m=\min(N_b,m)$ and replace $N_b$ with $N_m$ in Theorem \ref{th-1} (ii).
Since $N_m\uparrow N_b$ almost surely as $m\rightarrow\infty$, the result follows by the monotone convergence theorem.

(iii)

From (\ref{eq-23}), we have
\begin{align*}
\mathbb{E}(2N_b^2-5N_b+3) \leq (2n_1-3)(\mathbb{E}(N_b)-1) + 2v(\mathbb{E}(N_b)-1)\left(\sigma^2+\frac{(\mu_1-\mu_2)^2}{4}\right)
\end{align*}
implying
\begin{align*}
\mathbb{E}(N_b^2) \leq \left( n_1 + 1 + v\left(\sigma^2+\frac{(\mu_1-\mu_2)^2}{4}\right) \right)\mathbb{E}(N_b) -n_1 -v\left(\sigma^2+\frac{(\mu_1-\mu_2)^2}{4}\right).
\end{align*}

Applying Theorem \ref{th-1} (ii) again, the proof is completed.

\begin{flushright}
$\square$
\end{flushright}

\noindent\textbf{Proof of Theorem \ref{th-P(N<epsilona)}.}

The ideas of the proof are taken from \cite{Ghoshetal:1997}.
For convenience, we set $a=n_{req}=v\sigma^2$. We first consider the case where $\sigma \rightarrow \infty$. Choose $\sigma_1$ large enough such that for $\sigma\geq\sigma_1$, we have $\varepsilon a< [qa]-1$ for some $q$ in $(\varepsilon,1)$. Then, we have
\begin{align}\label{le-eq-1}
\mathbb{P}(N_b\leq \varepsilon a) &\leq \sum_{n=n_1}^{[\varepsilon a]}\mathbb{P}(N_b=n) \nonumber \\
&\leq \sum_{n=n_1}^{[qa]-1}\mathbb{P}(N_b=n) \nonumber \\
&\leq \sum_{n=n_1}^{[qa]-1}\mathbb{P}(n\geq v\hat{\sigma}^2_{n,blind}) \nonumber \\
&= \sum_{n=n_1}^{[qa]-1}\mathbb{P}\left(n\geq v\frac{\sigma^2}{2n-1}\chi_{2n-1}^2(\lambda)\right) \nonumber \\
&= \sum_{n=n_1}^{[qa]-1}\mathbb{P}\left(\chi_{2n-1}^2(\lambda)\leq \frac{n(2n-1)}{v\sigma^2}\right) \nonumber \\
&= \sum_{n=n_1}^{[qa]-1}\mathbb{P}\left(\chi_{2n-1}^2(\lambda)\leq \frac{n(2n-1)}{a}\right),
\end{align}
where $\chi_{2n-1}^2(\lambda)$ is noncentral chi-squared distribution with $2n-1$ degrees of freedom and noncentrality parameter $\lambda$.

From Markov's inequality and $h>0$, we have
\begin{align}\label{le-eq-2}
& \mathbb{P}\left(\chi_{2n-1}^2(\lambda)\leq \frac{n(2n-1)}{a}\right) \nonumber \\
&= \mathbb{P}\left(h\left(\frac{n(2n-1)}{a}-\chi_{2n-1}^2(\lambda)\right)\geq0\right) \nonumber \\
&= \mathbb{P}\left(\exp\left(h\left(\frac{n(2n-1)}{a}-\chi_{2n-1}^2(\lambda)\right)\right)\geq 1\right) \nonumber \\
&\leq \inf_{h>0}\mathbb{E}\left(\exp\left(h\left(\frac{n(2n-1)}{a}-\chi_{2n-1}^2(\lambda)\right)\right)\right) \ (\mbox{using moment generating function}) \nonumber \\
&= \inf_{h>0} \exp\left(h\left(\frac{n(2n-1)}{a}\right)\right)\exp\left(\frac{-\lambda h}{1+2h}\right)(1+2h)^{-(2n-1)/2} \nonumber \\
&\leq \inf_{h>0} \exp\left(h\left(\frac{n(2n-1)}{a}\right)\right)(1+2h)^{-(2n-1)/2} \nonumber \\
&= \exp\left(h_0\left(\frac{n(2n-1)}{a}\right)\right)(1+2h_0)^{-(2n-1)/2},
\end{align}
where
\begin{align*}
h_0 = \frac{1}{2}\left(\frac{a}{n}-1\right) \geq  \frac{1}{2}\left(\frac{1}{q}-1\right)>0
\end{align*}
for $0<q<1$ since $n\leq[qa]-1<qa-1<qa$. Now, from (\ref{le-eq-1}) and (\ref{le-eq-2}), we have
\begin{align}\label{le-eq-3}
\mathbb{P}(N_b\leq \varepsilon a)\leq \sum_{n=n_1}^{[qa]-1}\left(\exp\left(1-\frac{n}{a}\right)\frac{n}{a}\right)^{(2n-1)/2}.
\end{align}
Since $n\leq[qa]-1<qa-1<qa$, it follows that $\frac{n}{a}<q<1$. Note that the function $f(z)=z\exp(1-z)$ is increasing in $z$ for $0<z<1$. Now, let $\eta=\{q\exp(1-q)\}^{1/2}$. From (\ref{le-eq-3}), we have
\begin{align*}
\mathbb{P}(N_b\leq \varepsilon a) &\leq \sum_{n=n_1}^{[qa]-1}\eta^{(2n-1)-2(n_1-1)}\left(\exp\left(1-\frac{n}{a}\right)\frac{n}{a}\right)^{n_1-1} \\
&= a^{-(n_1-1)}\sum_{n=n_1}^{[qa]-1}\eta^{(2n-1)-2(n_1-1)}\left(\exp\left(1-\frac{n}{a}\right)n\right)^{n_1-1} \\
&\leq a^{-(n_1-1)}\sum_{n=n_1}^{[qa]-1}\eta^{(2n-1)-2(n_1-1)}(en)^{n_1-1} \\
&\leq a^{-(n_1-1)}\sum_{n=n_1}^{\infty}\eta^{(2n-1)-2(n_1-1)}(en)^{n_1-1} \\
&\leq Ka^{-(n_1-1)}
\end{align*}
for some constant $K(>0)$ not depending on $a$, since this sum can be shown to be convergent using the ratio rule of convergence.

Hence, for $n_1\geq2$, we have $\mathbb{P}(N_b\leq \varepsilon a) = \mathcal{O}(a^{-(n_1-1)})$ as $\sigma\rightarrow\infty$, since $a\rightarrow\infty$ as $\sigma\rightarrow\infty$.
The proof for $v \rightarrow \infty$ is similar.
This completes the proof.


\begin{flushright}
$\square$
\end{flushright}

\noindent\textbf{Proof of Theorem \ref{th-2}.}


(i)

Recall the stopping rule could be rewritten as
\begin{align*}
N_b &= \min\{ n=n_1,n_1+1,... | \hat{\sigma}^2_{n,blind}\leq n/v \} \\
  &= \min\{ n=n_1,n_1+1,... | \hat{\sigma}^2_{n,blind}/(\sigma^2+(\mu_1-\mu_2)^2/4) \leq n/(v(\sigma^2+(\mu_1-\mu_2)^2/4)) \}.
\end{align*}

From the definition of $\hat{\sigma}_{n,blind}^2$, we have
\begin{align*}
\frac{\hat{\sigma}_{n,blind}^2}{\sigma^2+(\mu_1-\mu_2)^2/4} \rightarrow 1 \quad almost \ surely \ as \ n\rightarrow\infty.
\end{align*}

Using this fact and the definition of $N_b$, we have $N_b\equiv N_b(\sigma)\rightarrow\infty$ almost surely as $\sigma\rightarrow\infty$.
By the monotone convergence theorem, it then follows that $\mathbb{E}(N_b)\rightarrow\infty$ as $\sigma\rightarrow\infty$. 
Readers are also referred to Lemma 1 in \cite{ChowRobbins:1965} for a related result.

(ii) 


Note that $N_b(\sigma)\rightarrow\infty$ almost surely as $\sigma\rightarrow\infty$.
Since we know
\begin{align*}
\frac{\hat{\sigma}_{n,blind}^2}{\sigma^2+(\mu_1-\mu_2)^2/4} \rightarrow 1 \quad almost \ surely \ as \ n\rightarrow\infty,
\end{align*}
it follows that
\begin{align*}
\hat{\sigma}^2_{N_b(\sigma),blind}/\sigma^2  \rightarrow 1 \quad almost \ surely \  as \ \sigma\rightarrow\infty,
\end{align*}
and
\begin{align*}
\hat{\sigma}^2_{N_b(\sigma)-1,blind}/\sigma^2  \rightarrow 1 \quad almost \ surely \ as \ \sigma\rightarrow\infty
\end{align*}
by Theorem 2.1 of \cite{Gut:2009}. 

We now use the inequality
\begin{align*}
v\hat{\sigma}^2_{N_b,blind} \leq N_b \leq n_1 + v\hat{\sigma}^2_{N_b-1,blind} \quad almost \ surely .
\end{align*}
Dividing both sides by $n_{req}=v\sigma^2$ and making $\sigma\rightarrow\infty$, we complete the proof.

(iii)

Since Theorem \ref{th-2} (ii) and Fatou's lemma, it follows that
\begin{align*}
\liminf_{\sigma\rightarrow\infty} \mathbb{E}\left(\frac{N_b}{n_{req}}\right) \geq \mathbb{E}\left(\liminf_{\sigma\rightarrow\infty} \frac{N_b}{n_{req}}\right) = 1.
\end{align*}
From Theorem \ref{th-1} (ii), we have
\begin{align*}
\mathbb{E}\left(\frac{N_b}{n_{req}}\right) \leq \frac{n_1}{v\sigma^2}+1+\frac{(\mu_1-\mu_2)^2}{4\sigma^2}
\end{align*}
which leads to
\begin{align*}
\limsup_{\sigma\rightarrow\infty} \mathbb{E}\left(\frac{N_b}{n_{req}}\right) \leq 1.
\end{align*}

Thus, we have
\begin{align*}
1 \leq \liminf_{\sigma\rightarrow\infty} \mathbb{E}\left(\frac{N_b}{n_{req}}\right) \leq \limsup_{\sigma\rightarrow\infty} \mathbb{E}\left(\frac{N_b}{n_{req}}\right) \leq 1,
\end{align*}
which completes the proof.

(iv)

Since Theorem \ref{th-2} (ii) and Fatou's lemma, it follows that
\begin{align*}
\liminf_{\sigma\rightarrow\infty} \mathbb{E}\left(\left(\frac{N_b}{n_{req}}\right)^2\right) \geq \mathbb{E}\left(\liminf_{\sigma\rightarrow\infty} \left(\frac{N_b}{n_{req}}\right)^2\right) = 1.
\end{align*}
From Theorem \ref{th-1} (iii), we have
\begin{align*}
\mathbb{E}\left(\left(\frac{N_b}{n_{req}}\right)^2\right) \leq \left( \frac{n_1}{v\sigma^2}+1+\frac{(\mu_1-\mu_2)^2}{4\sigma^2}\right)^2
\end{align*}
which leads to
\begin{align*}
\limsup_{\sigma\rightarrow\infty} \mathbb{E}\left(\left(\frac{N_b}{n_{req}}\right)^2\right) \leq 1.
\end{align*}

Thus, we have
\begin{align*}
1 \leq \liminf_{\sigma\rightarrow\infty} \mathbb{E}\left(\left(\frac{N_b}{n_{req}}\right)^2\right) \leq \limsup_{\sigma\rightarrow\infty} \mathbb{E}\left(\left(\frac{N_b}{n_{req}}\right)^2\right) \leq 1,
\end{align*}
which completes the proof.

(v)

Since we have the basic inequality
\begin{align*}
v\hat{\sigma}^2_{N_b,blind} \leq N_b \leq n_1 + v\hat{\sigma}^2_{N_b-1,blind} \quad almost \  surely,
\end{align*}
it follows that
\begin{align}\label{eq-333sigma}
\sqrt{n_{req}}\left( \frac{\hat{\sigma}^2_{N_b,blind}}{\sigma^2}-1 \right) \leq \frac{N_b-n_{req}}{\sqrt{n_{req}}} \leq \frac{n_1}{\sqrt{n_{req}}} + \sqrt{n_{req}}\left( \frac{\hat{\sigma}^2_{N_b-1,blind}}{\sigma^2}-1 \right),
\end{align}
where $n_{req}=v\sigma^2$.

Recall that $\hat{\sigma}^2_{n,blind}/\sigma^2$ follows $\chi^2_{2n-1}(\lambda)/(2n-1)$ where $\chi^2_{2n-1}(\lambda)$ is noncentral chi-squared distribution with $2n-1$ degrees of freedom and noncentrality parameter $\lambda=n(\mu_1-\mu_2)^2/(2\sigma^2)$.

Note that the left hand side of (\ref{eq-333sigma}) can be rewritten as
\begin{align*}
& \sqrt{n_{req}}\left(\frac{\hat{\sigma}^2_{N_b,blind}}{\sigma^2}-1\right) \\
&= \sqrt{n_{req}}\left(\frac{\chi^2_{2N_b-1}(\lambda)}{2N_b-1}-1\right) \\
&= \sqrt{n_{req}}\left(\frac{\chi^2_{2N_b-1}(\lambda)-(2N_b-1+\lambda)}{2N_b-1}+\frac{\lambda}{2N_b-1}\right) \\
&= \sqrt{n_{req}}\left(\frac{\sqrt{2(2N_b-1+2\lambda)}}{2N_b-1}\frac{\chi^2_{2N_b-1}(\lambda)-(2N_b-1+\lambda)}{\sqrt{2(2N_b-1+2\lambda)}}+\frac{\lambda}{2N_b-1}\right) \\
&= \sqrt{\frac{n_{req}}{N_b}}\frac{\sqrt{2N_b(2N_b-1+2\lambda)}}{2N_b-1}\frac{\chi^2_{2N_b-1}(\lambda)-(2N_b-1+\lambda)}{\sqrt{2(2N_b-1+2\lambda)}}+\frac{N_b}{2N_b-1}\sqrt{n_{req}}\frac{(\mu_1-\mu_2)^2}{2\sigma^2} \\
&= \sqrt{\frac{n_{req}}{N_b}}\sqrt{\frac{(4+2(\mu_1-\mu_2)^2/\sigma^2)N_b^2-2N_b}{4N_b^2-4N_b+1}}\frac{\chi^2_{2N_b-1}(\lambda)-(2N_b-1+\lambda)}{\sqrt{2(2N_b-1+2\lambda)}}+\frac{N_b}{2N_b-1}\sqrt{n_{req}}\frac{(\mu_1-\mu_2)^2}{2\sigma^2} \\
&= \sqrt{\frac{n_{req}}{N_b}}\sqrt{\frac{4+2(\mu_1-\mu_2)^2/\sigma^2-2/N_b}{4-4/N_b+1/N_b^2}}\frac{\chi^2_{2N_b-1}(\lambda)-(2N_b-1+\lambda)}{\sqrt{2(2N_b-1+2\lambda)}}+\frac{N_b}{2N_b-1}\sqrt{v\sigma^2}\frac{(\mu_1-\mu_2)^2}{2\sigma^2}.
\end{align*}

Thus, using Anscombe's Theorem (see Theorem 2.7.1 in \cite{Ghoshetal:1997}) and Theorem \ref{th-2} (ii), we have
\begin{align*}
\sqrt{n_{req}}\left(\frac{\hat{\sigma}^2_{N_b,blind}}{\sigma^2}-1\right) \overset{d}{\rightarrow} N(0,1) \quad \mbox{and} \quad \sqrt{n_{req}}\left(\frac{\hat{\sigma}^2_{N_b-1,blind}}{\sigma^2}-1\right) \overset{d}{\rightarrow} N(0,1)
\end{align*}
as $\sigma\rightarrow\infty$. The result follows now from (\ref{eq-333sigma}).


\begin{flushright}
$\square$
\end{flushright}

\noindent\textbf{Proof of Theorem \ref{th-3}.}

(i)

From the definition of $\hat{\sigma}_{n,blind}^2$, we have
\begin{align*}
\hat{\sigma}_{n,blind}^2 \rightarrow \sigma^2+(\mu_1-\mu_2)^2/4 \quad almost \ surely \ as \ n\rightarrow\infty.
\end{align*}

Using this fact and the definition of $N_b$, we have $N_b\equiv N_b(v)\rightarrow\infty$ almost surely as $v\rightarrow\infty$.
From the monotone convergence theorem, we have $\mathbb{E}(N_b)\rightarrow\infty$ as $v\rightarrow\infty$.

(ii)

Note that $N_b(v)\rightarrow\infty$ almost surely as $v\rightarrow\infty$.
Hence, by Theorem 2.1 of \cite{Gut:2009}, we have
\begin{align*}
\hat{\sigma}^2_{N_b(v),blind}  \rightarrow \sigma^2 + (\mu_1-\mu_2)^2/4 \quad almost \ surely \ as \ v\rightarrow\infty,
\end{align*}
and
\begin{align*}
\hat{\sigma}^2_{N_b(v)-1,blind}  \rightarrow \sigma^2 + (\mu_1-\mu_2)^2/4 \quad almost \ surely \ as \ v\rightarrow\infty.
\end{align*}

We now use the inequality
\begin{align*}
v\hat{\sigma}^2_{N_b,blind} \leq N_b \leq n_1 + v\hat{\sigma}^2_{N-1,blind} \quad almost \ surely
\end{align*}
Dividing both sides by $n_{req}=v\sigma^2$ and making $v\rightarrow\infty$, we complete the proof.

(iii)

Since Theorem \ref{th-3} (ii) and Fatou's lemma, it follows that
\begin{align*}
\liminf_{v\rightarrow\infty} \mathbb{E}\left(\frac{N_b}{n_{req}}\right) \geq \mathbb{E}\left(\liminf_{v\rightarrow\infty} \frac{N_b}{n_{req}}\right) = 1+\frac{(\mu_1-\mu_2)^2}{4\sigma^2}.
\end{align*}
From Theorem \ref{th-1} (ii), we have
\begin{align*}
\mathbb{E}\left(\frac{N_b}{n_{req}}\right) \leq \frac{n_1}{v\sigma^2} + 1+\frac{(\mu_1-\mu_2)^2}{4\sigma^2}
\end{align*}
which leads to
\begin{align*}
\limsup_{v\rightarrow\infty} \mathbb{E}\left(\frac{N_b}{n_{req}}\right) \leq 1+\frac{(\mu_1-\mu_2)^2}{4\sigma^2}.
\end{align*}

Thus, we have
\begin{align*}
1+\frac{(\mu_1-\mu_2)^2}{4\sigma^2} \leq \liminf_{v\rightarrow\infty} \mathbb{E}\left(\frac{N_b}{n_{req}}\right) \leq \limsup_{v\rightarrow\infty} \mathbb{E}\left(\frac{N_b}{n_{req}}\right) \leq 1+\frac{(\mu_1-\mu_2)^2}{4\sigma^2},
\end{align*}
which completes the proof.

(iv)

Since Theorem \ref{th-3} (ii) and Fatou's lemma, it follows that
\begin{align}
\liminf_{v\rightarrow\infty} \mathbb{E}\left(\left(\frac{N_b}{n_{req}}\right)^2\right) \geq \mathbb{E}\left(\liminf_{v\rightarrow\infty} \left(\frac{N}{n_{req}}\right)^2\right) = \left(1+\frac{(\mu_1-\mu_2)^2}{4\sigma^2}\right)^2.
\end{align}
From Theorem \ref{th-1} (iii), we have
\begin{align}
\mathbb{E}\left(\left(\frac{N_b}{n_{req}}\right)^2\right) \leq \left( \frac{n_1}{v\sigma^2} + 1+\frac{(\mu_1-\mu_2)^2}{4\sigma^2} \right)^2
\end{align}
which leads to
\begin{align}
\limsup_{v\rightarrow\infty} \mathbb{E}\left(\left(\frac{N_b}{n_{req}}\right)^2\right) \leq \left(1+\frac{(\mu_1-\mu_2)^2}{4\sigma^2}\right)^2.
\end{align}

Thus, we have
\begin{align}
\left(1+\frac{(\mu_1-\mu_2)^2}{4\sigma^2}\right)^2 \leq \liminf_{v\rightarrow\infty} \mathbb{E}\left(\left(\frac{N_b}{n_{req}}\right)^2\right) \leq \limsup_{v\rightarrow\infty} \mathbb{E}\left(\left(\frac{N_b}{n_{req}}\right)^2\right) \leq \left(1+\frac{(\mu_1-\mu_2)^2}{4\sigma^2}\right)^2,
\end{align}
which completes the proof.

(v)

Since we have the basic inequality
\begin{align*}
v\hat{\sigma}^2_{N_b,blind} \leq N_b \leq n_1 + v\hat{\sigma}^2_{N_b-1,blind} \quad almost \ surely,
\end{align*}
it follows that
\begin{align}\label{eq-333}
\sqrt{n_{req}}\left(\frac{\hat{\sigma}^2_{N_b,blind}}{\sigma^2}-1-\frac{(\mu_1-\mu_2)^2}{4\sigma^2}\right) &\leq \frac{N_b-n_{req}(1+(\mu_1-\mu_2)^2/4\sigma^2)}{\sqrt{n_{req}}} \nonumber \\
&\leq \frac{n_1}{\sqrt{n_{req}}} + \sqrt{n_{req}}\left(\frac{\hat{\sigma}^2_{N_b-1,blind}}{\sigma^2}-1-\frac{(\mu_1-\mu_2)^2}{4\sigma^2}\right),
\end{align}
where $n_{req}=v\sigma^2$

Recall that $\hat{\sigma}^2_{n,blind}/\sigma^2$ follows $\chi^2_{2n-1}(\lambda)/(2n-1)$ where $\chi^2_{2n-1}(\lambda)$ is noncentral chi-squared distribution with degree of freedom $2n-1$ and the noncentrality parameter $\lambda=n(\mu_1-\mu_2)^2/(2\sigma^2)$.

Note that the left hand side of (\ref{eq-333}) can be rewritten as
\begin{align*}
& \sqrt{n_{req}}\left(\frac{\hat{\sigma}^2_{N_b,blind}}{\sigma^2}-1-\frac{(\mu_1-\mu_2)^2}{4\sigma^2}\right) \\
&= \sqrt{n_{req}}\left(\frac{\chi^2_{2N_b-1}(\lambda)}{2N_b-1}-1-\frac{(\mu_1-\mu_2)^2}{4\sigma^2}\right) \\
&= \sqrt{n_{req}}\left(\frac{\chi^2_{2N_b-1}(\lambda)-(2N_b-1+\lambda)}{2N_b-1}+\frac{1}{2N_b-1}\frac{(\mu_1-\mu_2)^2}{4\sigma^2}\right) \\
&= \sqrt{n_{req}}\left(\frac{\sqrt{2(2N_b-1+2\lambda)}}{2N_b-1}\frac{\chi^2_{2N_b-1}(\lambda)-(2N_b-1+\lambda)}{\sqrt{2(2N_b-1+2\lambda)}}+\frac{1}{2N_b-1}\frac{(\mu_1-\mu_2)^2}{4\sigma^2}\right) \\
&= \sqrt{\frac{n_{req}}{N_b}}\frac{\sqrt{2N_b(2N_b-1+2\lambda)}}{2N_b-1}\frac{\chi^2_{2N_b-1}(\lambda)-(2N_b-1+\lambda)}{\sqrt{2(2N_b-1+2\lambda)}}+\frac{\sqrt{n_{req}}}{2N_b-1}\frac{(\mu_1-\mu_2)^2}{4\sigma^2} \\
&= \sqrt{\frac{n_{req}}{N_b}}\sqrt{\frac{(4+2(\mu_1-\mu_2)^2/\sigma^2)N_b^2-2N_b}{4N_b^2-4N_b+1}}\frac{\chi^2_{2N_b-1}(\lambda)-(2N_b-1+\lambda)}{\sqrt{2(2N_b-1+2\lambda)}}+\frac{\sqrt{n_{req}}}{2N_b-1}\frac{(\mu_1-\mu_2)^2}{4\sigma^2} \\
&= \sqrt{\frac{n_{req}}{N_b}}\sqrt{\frac{4+2(\mu_1-\mu_2)^2/\sigma^2-2/N_b}{4-4/N_b+1/N_b^2}}\frac{\chi^2_{2N_b-1}(\lambda)-(2N_b-1+\lambda)}{\sqrt{2(2N_b-1+2\lambda)}}+\frac{\sqrt{v\sigma^2}}{2N_b-1}\frac{(\mu_1-\mu_2)^2}{4\sigma^2}.
\end{align*}

Thus, using Anscombe's Theorem (see Theorem 2.7.1 in \cite{Ghoshetal:1997}) and Theorem \ref{th-3} (ii), we have
\begin{align*}
\sqrt{n_{req}}\left(\frac{\hat{\sigma}^2_{N_b,blind}}{\sigma^2}-1-\frac{(\mu_1-\mu_2)^2}{4\sigma^2}\right) \overset{d}{\rightarrow} N(0,\frac{4\sigma^2+2(\mu_1-\mu_2)^2}{4\sigma^2+(\mu_1-\mu_2)^2})
\end{align*}
and
\begin{align*}
\sqrt{n_{req}}\left(\frac{\hat{\sigma}^2_{N_b-1,blind}}{\sigma^2}-1-\frac{(\mu_1-\mu_2)^2}{4\sigma^2}\right) \overset{d}{\rightarrow} N(0,\frac{4\sigma^2+2(\mu_1-\mu_2)^2}{4\sigma^2+(\mu_1-\mu_2)^2})
\end{align*}
as $v\rightarrow\infty$. This result follows now from (\ref{eq-333}).

\begin{flushright}
$\square$
\end{flushright}


\end{document}